\documentclass{article}

\usepackage{arxiv}

\usepackage[utf8]{inputenc} 
\usepackage[T1]{fontenc}    

\usepackage{hyperref}       
\usepackage{url}            
\usepackage{booktabs}       
\usepackage{amsfonts}       
\usepackage{amsmath}
\usepackage{amssymb}
\usepackage{microtype}      
\usepackage{graphicx}
\usepackage{lipsum}
\usepackage[numbers]{natbib} 
\usepackage{doi}
\usepackage{float}
\usepackage[section]{placeins}
\FloatBarrier

\usepackage{stmaryrd}
\hypersetup{
  colorlinks = true,
  linkcolor  = black,  
  citecolor  = blue,   
  urlcolor   = blue    
}

\title{Adaptive Wavelet–Galerkin Modelling of Heat Conduction in Heterogeneous Composite Materials}

\author{ \href{https://orcid.org/0000-0001-7402-4468}{\includegraphics[scale=0.06]{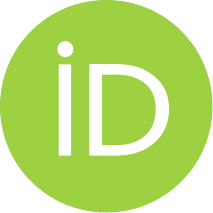}\hspace{1mm}Taylan Demir}\thanks{Use footnote for providing further
		information about author (webpage, alternative
		address)---\emph{not} for acknowledging funding agencies.} \\
	Department of Mathematics\\
	  Ankara University\\
	Ankara, Turkey \\
	\texttt{taylandemir@ankara.edu.tr} \\
	\And
	\href{https://orcid.org/0009-0005-6060-0547}{\includegraphics[scale=0.06]{orcid.pdf}\hspace{1mm}Atakan Koçyiğit} \\
	Department of Electrical and Electronics Engineering\\
	  Atılım University\\
	Ankara, Turkey \\
	\texttt{atakan.kocyigit@megetek.com.tr} \\
}

\renewcommand{\shorttitle}{\textit{arXiv} Template}

\hypersetup{
pdftitle={A template for the arxiv style},
pdfsubject={q-bio.NC, q-bio.QM},
pdfauthor={David S.~Hippocampus, Elias D.~Striatum},
pdfkeywords={First keyword, Second keyword, More},
}

\begin{document}
\maketitle
\begin{abstract}
We present an adaptive wavelet Galerkin method for transient heat conduction in heterogeneous composite materials. The approach combines multiresolution wavelet bases with an implicit time discretization to efficiently resolve sharp temperature gradients near material interfaces and boundary layers. Adaptive refinement is driven by wavelet coefficients, significantly reducing the number of degrees of freedom compared to uniform discretizations. Numerical examples demonstrate accurate resolution of layered, inclusion-based, and functionally graded composites with improved computational efficiency.
\end{abstract}

\keywords{Wavelet--Galerkin method, heat conduction, composite materials,
multiresolution analysis, adaptive numerical methods, thermal analysis}
\section{Introduction}
Engineered composite materials consist of multiple phases at various length scales that come together to create characteristics including stiffness, strength, weight, and thermal performance that are impossible to achieve with conventional metals or polymeric materials [1,2]. There are numerous applications for composite materials; for example, in the fields of aerospace structures, automotive components, electronic packaging and energy systems, thermal conductivity (efficiency of transfer and redistribution of heat) is as significant for design purposes as it is for mechanical function. These types of applications will typically cause localised heating and resultant degradation of matrix due to induced residual stresses, which in turn can accelerate fatigue or delamination. Reliable prediction of temperature fields is therefore important for safe design, durability assessment and for optimising composite structures. The classical heat equation is currently used at the continuum scale for modeling heat transport within composite structures; however, in nature there are many instances where heat transport does not occur evenly throughout the material because of spatially varying and generally anisotropic thermal conductivities. Each of the individual plies, inclusions, and matrix phases will have different thermal conductivity tensors, so that there will be large contrasts between these components at their point of interface and large non-uniform temperature gradients will exist across the thickness of the composite. Even if there is considerable time and expense involved with experimental thermal characterization, experiments can only be conducted under a limited set of load and boundary conditions [3]. Consequently, there is a continuous requirement for numerical methods that are able to represent heat diffusion in complicated composite shapes through adequate accuracy and with a computational expense that is compatible with parametric study and design loop requirements. There are three primary types of discretisation methods that are commonly used to solve static heat conduction problems in multi-layered and fibre reinforced composite materials: Finite-difference, finite-volume and finite element. On refined mesh sizes, these three methods can provide accurate representations of thin films, layered materials and/ or heat generation points with high gradients of temperature [4]. In examples of highly heterogeneous or complex geometries, however, the mesh is required to be aligned with the heterogeneities locally across the interface, as well as within every region of large change in the heat conduction behaviour of the material; therefore, the solution will typically require an enormous algebraic system of equations to be solved, and as a result of geometry and/or change in material distribution from the prior step, the need to re-mesh repeatedly will exist. In addition, the need to reiterate through successive changes to achieve the proper balance of computationally efficient versus accurate will place significant demands on the operator to perform manual tuning. This problem would be exacerbated when trying to understand several different types of behaviour, from a sharp jump in the thermal conductivity of materials (layered or particulate composites) to the smooth progressive change in thermal properties, all within one set of numerical analysis programme. Alternatively, the use of wavelet basis functions will provide a method of representing the spatial discretisation of the partial differential equations with inherent multiscale structures. Wavelet basis functions that have finite support, such as those developed by Daubechies, have properties of concurrent localisation in both scale and space and provide, when using a Galerkin approach, either sparse or near sparse representations of the differentiation operators [5].
Wavelet-based approaches have been effectively used for many heat conduction problems, including inverse problems and sideways problems, by utilizing the multiresolution representation both for regularisation and for providing an intuitive method for adaptively refining the solution space over time [6,7]. These works show that wavelet-Galerkin methods permit the accurate reproduction of sharp spatial feature resolution and control over the impact of high-frequency noise on the solution. The use of adaptive wavelet–Galerkin methods to model heat conduction in strongly variable composites is not well developed and, as compared to the extensive body of literature on finite–element techniques for modelling heat conduction, there have been only a few papers published. The majority of work has focused on one-dimensional cases or inverse problems without creating a unified treatment that describes all of the possible (i) piecewise constant thermal conductivities with sharp boundaries, (ii) curved particles present within a blend of materials, and (iii) graded mixtures of materials that exhibit temperature gradients. However, these three general types of thermal conductivities are present in numerous applications in engineering. Composite materials such as coated and layered structures fall into the first category, reinforced composites made up of particles or fibres can generally be categorized under the second category, and functionally graded materials fall into the third category. It is therefore necessary to create a comprehensive numerical technique that can simultaneously model all three types of toxic waste management practices and efficiently capture the temperature changes at the interface where the temperature gradient changes occur. This study's objective is to create and test the application of the adaptive wavelet-Galerkin method to numerically solve the heat equation in heterogeneous composite materials based on practical design applications. This will be accomplished by utilizing a multi-resolution basis on a fixed computational domain to create a wavelet basis set and utilizing the size of the wavelet coefficients as an indicator of the error to adaptively refine the mesh in the areas where the temperature has high gradients, e.g., at the interfaces of materials, near the circular inclusions that are embedded, and across regions with high thermal gradients. In addition, the study will show the performance of the method on three representative examples of layered composites, a matrix with a circular inclusion, and a functionally graded composite so that the adaptive refinement pattern is representative of the micro-structure of the materials and can be utilized for thermal design studies.
\section{Mathematical Model}

We consider a bounded Lipschitz domain $\Omega \subset \mathbb{R}^d$, $d \in \{2,3\}$, occupied by a heterogeneous composite material, over a finite time interval $(0,t_{\mathrm{f}}]$. 
The unknown temperature field is denoted by $T(x,t)$, where $x \in \Omega$ and $t \in (0,t_{\mathrm{f}}]$. 
At the macroscopic level, heat transport is described by the transient anisotropic heat equation
\begin{equation}
  \rho(x)c_{p}(x)\,\frac{\partial T}{\partial t}(x,t)
  - \nabla\!\cdot\!\bigl(\mathbf{K}(x)\nabla T(x,t)\bigr)
  = f(x,t)
  \qquad \text{in } \Omega \times (0,t_{\mathrm{f}}],
  \label{eq:heat-strong}
\end{equation}
where $\rho(x)$ is the mass density, $c_{p}(x)$ is the specific heat capacity, $\mathbf{K}(x)$ is the (possibly anisotropic) thermal conductivity tensor and $f(x,t)$ denotes a prescribed volumetric heat source [3,8]. 
In composite materials these coefficients may vary strongly in space and exhibit jumps across material interfaces.

The boundary $\partial\Omega$ is decomposed into three pairwise disjoint parts,
\[
  \partial\Omega = \Gamma_{D} \,\dot{\cup}\, \Gamma_{N} \,\dot{\cup}\, \Gamma_{R},
\]
on which we impose, respectively, Dirichlet, Neumann and Robin boundary conditions. 
On $\Gamma_{D}$ the temperature is prescribed as
\begin{equation}
  T(x,t) = g_{D}(x,t)
  \qquad \text{for } x \in \Gamma_{D}, \; t \in (0,t_{\mathrm{f}}],
  \label{eq:bc-dirichlet}
\end{equation}
while on $\Gamma_{N}$ a given normal heat flux $g_{N}$ is specified,
\begin{equation}
  -\mathbf{n}(x)\cdot\mathbf{K}(x)\nabla T(x,t) = g_{N}(x,t)
  \qquad \text{for } x \in \Gamma_{N}, \; t \in (0,t_{\mathrm{f}}],
  \label{eq:bc-neumann}
\end{equation}
with $\mathbf{n}(x)$ denoting the outward unit normal to $\partial\Omega$. 
On the remaining part of the boundary, $\Gamma_{R}$, convective exchange with an ambient medium of temperature $T_{\infty}(x,t)$ is modelled by
\begin{equation}
  -\mathbf{n}(x)\cdot\mathbf{K}(x)\nabla T(x,t)
  = h(x)\,\bigl(T(x,t)-T_{\infty}(x,t)\bigr)
  \qquad \text{for } x \in \Gamma_{R}, \; t \in (0,t_{\mathrm{f}}],
  \label{eq:bc-robin}
\end{equation}
where $h(x) \ge 0$ denotes the heat transfer coefficient on $\Gamma_{R}$ [3]. 
The initial state of the body is prescribed by
\begin{equation}
  T(x,0) = T_{0}(x)
  \qquad \text{for } x \in \Omega,
  \label{eq:initial}
\end{equation}
with $T_{0}$ a given temperature distribution.

To represent the heterogeneous character of the composite, we assume that $\Omega$ is partitioned into subdomains corresponding to distinct material phases,
\begin{equation}
  \Omega = \Omega_{\mathrm{m}} \,\dot{\cup}\,
           \Bigl(\bigcup_{k=1}^{N_{\mathrm{inc}}} \Omega_{\mathrm{inc}}^{(k)}\Bigr)
           \,\dot{\cup}\, \Omega_{\mathrm{fgm}},
  \label{eq:domain-decomposition}
\end{equation}
where $\Omega_{\mathrm{m}}$ denotes the matrix region, $\Omega_{\mathrm{inc}}^{(k)}$ are inclusion domains (for instance circular or elliptical particles) and $\Omega_{\mathrm{fgm}}$ is an optional functionally graded region. 
In the matrix and in each inclusion, the material parameters are taken to be piecewise constant,
\[
  \rho(x) = \rho_{j}, 
  \quad
  c_{p}(x) = c_{p,j},
  \quad
  \mathbf{K}(x) = \mathbf{K}_{j}
  \qquad \text{for } x \in \Omega_{j},
\]
with $\Omega_{j} \in \bigl\{\Omega_{\mathrm{m}},\Omega_{\mathrm{inc}}^{(1)},\dots,\Omega_{\mathrm{inc}}^{(N_{\mathrm{inc}})}\bigr\}$, 
while in the graded region $\Omega_{\mathrm{fgm}}$ the conductivity tensor may vary smoothly in space according to a prescribed grading law, for example along a coating thickness or through the plate thickness direction [9]. 
Across each internal interface between two phases we impose perfect thermal contact, which at the macroscopic level amounts to continuity of temperature and of the normal component of the heat flux,
\begin{equation}
  \llbracket T \rrbracket_{\Gamma_{\mathrm{int}}} = 0,
  \qquad
  \llbracket \mathbf{n}\cdot\mathbf{K}\nabla T \rrbracket_{\Gamma_{\mathrm{int}}} = 0,
  \label{eq:interface}
\end{equation}
for all internal interfaces $\Gamma_{\mathrm{int}}$ separating adjacent subdomains. 
Here $\llbracket \cdot \rrbracket_{\Gamma_{\mathrm{int}}}$ denotes the jump across the interface and $\mathbf{n}$ is chosen consistently with the outward normal of one of the subdomains.

From a mathematical standpoint we assume that, for almost every $x \in \Omega$, the conductivity tensor $\mathbf{K}(x)$ is symmetric and uniformly positive definite, i.e.
\begin{equation}
  k_{\min} \, |\xi|^{2}
  \le \xi^{\mathsf{T}} \mathbf{K}(x)\,\xi
  \le k_{\max} \, |\xi|^{2}
  \qquad \forall\, \xi \in \mathbb{R}^{d},
  \label{eq:uniform-ellipticity}
\end{equation}
for some constants $0 < k_{\min} \le k_{\max} < \infty$. 
Together with standard regularity assumptions on the data $f$, $g_{D}$, $g_{N}$, $T_{\infty}$ and $T_{0}$, these conditions guarantee the well-posedness of the initial–boundary value problem \eqref{eq:heat-strong}--\eqref{eq:initial} in appropriate Sobolev spaces [10]. 
This strong formulation will serve as the starting point for the wavelet--Galerkin spatial discretisation developed in the next section.
\section{Wavelet--Galerkin Method}
\label{sec:wavelet-galerkin}
An adaptive wavelet-Galerkin method is presented here as a method for discretising (i.e. approximating) problem (1)-(5). We outline the suitable Sobolev space in which to formulate weakly the heat equation, and substitute the standard Hilbert basis with that of an adaptive multiresolution wavelet basis keyed to the boundary conditions. The system now becomes a set of semi-discrete ODEs represented by a sparse (or compressible) stiffness matrix. The adaptive nature of this scheme is apparent in the way that the wavelet coefficients are dynamically chosen to achieve higher resolution at the interfaces, inclusions and thermal boundary layers of the domain.

\subsection{Weak formulation}

Let
\[
  V := \bigl\{ v \in H^{1}(\Omega) \; : \; v = 0 \text{ on } \Gamma_{D} \bigr\},
\]
equipped with the standard $H^{1}$-norm, and define the bilinear form
$a : V \times V \to \mathbb{R}$ and the (time-dependent) linear functional
$\ell(\cdot,t) : V \to \mathbb{R}$ by
\begin{align}
  a(u,v)
  &:= \int_{\Omega} \nabla v(x)^{\mathsf{T}} \mathbf{K}(x)\nabla u(x)\,\mathrm{d}x
      + \int_{\Gamma_{R}} h(x)\,u(x)v(x)\,\mathrm{d}s(x),
  \label{eq:bilinear-form}
  \\
  \ell(v,t)
  &:= \int_{\Omega} f(x,t)\,v(x)\,\mathrm{d}x
      + \int_{\Gamma_{N}} g_{N}(x,t)\,v(x)\,\mathrm{d}s(x)
      + \int_{\Gamma_{R}} h(x)\,T_{\infty}(x,t)\,v(x)\,\mathrm{d}s(x),
  \label{eq:functional}
\end{align}
for $u,v \in V$ and $t \in (0,t_{\mathrm{f}}]$. 
Moreover, we introduce the weighted mass bilinear form
\begin{equation}
  m(u,v) := \int_{\Omega} \rho(x)c_{p}(x)\,u(x)v(x)\,\mathrm{d}x,
  \qquad u,v \in V.
  \label{eq:mass-form}
\end{equation}
Under the assumptions stated in Section 2, in
particular the uniform ellipticity condition \eqref{eq:uniform-ellipticity},
the form $a(\cdot,\cdot)$ is continuous and coercive on $V$ and
$m(\cdot,\cdot)$ is continuous and symmetric positive definite [10].

The weak formulation of the initial--boundary value problem reads as
follows: find
\[
  T(\cdot,\cdot) \in L^{2}\bigl(0,t_{\mathrm{f}};V\bigr)
  \quad \text{with} \quad
  \partial_{t}T(\cdot,\cdot) \in L^{2}\bigl(0,t_{\mathrm{f}};V^{*}\bigr)
\]
such that
\begin{align}
  m\bigl(\partial_{t}T(\cdot,t), v\bigr) + a\bigl(T(\cdot,t),v\bigr)
  &= \ell(v,t)
  \qquad \forall v \in V,\; t \in (0,t_{\mathrm{f}}],
  \label{eq:weak-form}
  \\
  T(\cdot,0) &= T_{0} \quad \text{in } \Omega.
  \label{eq:weak-initial}
\end{align}
This functional setting will serve as the starting point for the
wavelet--Galerkin discretisation discussed below.

\subsection{Multiresolution wavelet bases}

To construct a wavelet basis for $V$ we start from a nested sequence of
finite-dimensional subspaces
\[
  V_{0} \subset V_{1} \subset \dots \subset V_{J} \subset V \subset H^{1}(\Omega),
\]
where each $V_{j}$ is spanned by suitable scaling functions supported at a
characteristic length scale $2^{-j}$ [11,12]. 
The index $J$ denotes the finest resolution level admitted in the
computation and is chosen large enough so that $V_{J}$ can approximate the
exact solution with the desired accuracy. 
Associated with this multiresolution structure is a family of wavelet
spaces
\[
  W_{j} \subset L^{2}(\Omega),
  \qquad j = 0,1,\dots,J-1,
\]
such that
\[
  V_{j+1} = V_{j} \oplus W_{j},
  \qquad j=0,\dots,J-1,
\]
and, in the limit $J\to\infty$,
\[
  \overline{\bigcup_{j\ge 0} V_{j}}^{\,L^{2}(\Omega)} = L^{2}(\Omega).
\]
In practice we employ compactly supported, orthogonal (or biorthogonal)
wavelets of Daubechies type, modified near the boundary so that the
resulting basis functions satisfy the homogeneous Dirichlet conditions on
$\Gamma_{D}$ [12,13]. 
The boundary adaptation ensures that the discrete trial space
inherits the correct essential boundary conditions without the need
for explicit constraints.

Let $\{\psi_{\lambda} : \lambda \in \Lambda_{J}\}$ denote the collection
of all wavelets up to level $J-1$ together with the coarsest-level scaling
functions in $V_{0}$, where $\Lambda_{J}$ is a suitable multi-index set
encoding position, orientation and scale. 
It is well known that, under mild conditions, this system forms a
Riesz basis of $L^{2}(\Omega)$ (and of appropriate Sobolev spaces), so
that any $v \in V$ can be written as a convergent expansion
\begin{equation}
  v(x) = \sum_{\lambda \in \Lambda_{J}} v_{\lambda}\,\psi_{\lambda}(x),
  \qquad x \in \Omega,
  \label{eq:wavelet-expansion}
\end{equation}
with square-summable coefficients $(v_{\lambda})_{\lambda \in \Lambda_{J}}$ [11,13]. 
The localisation properties of the functions $\psi_{\lambda}$ entail that
only a small number of basis functions have support in a given region of
$\Omega$, which is the key to the sparsity of the resulting stiffness
matrix.

\subsection{Wavelet--Galerkin semi-discretisation}

We now approximate the temperature field $T(\cdot,t)$ by a finite
wavelet expansion of the form
\begin{equation}
  T_{J}(x,t) = \sum_{\lambda \in \Lambda_{J}} u_{\lambda}(t)\,\psi_{\lambda}(x),
  \qquad x \in \Omega,\; t \in (0,t_{\mathrm{f}}],
  \label{eq:temp-expansion}
\end{equation}
where the functions $u_{\lambda}(t)$ are unknown time-dependent
coefficients. 
Inserting \eqref{eq:temp-expansion} into the weak formulation
\eqref{eq:weak-form} and taking test functions $v$ from the span of
$\{\psi_{\lambda} : \lambda \in \Lambda_{J}\}$ leads to a system of
ordinary differential equations,
\begin{equation}
  \mathbf{M}\,\dot{\mathbf{u}}(t) + \mathbf{K}\,\mathbf{u}(t)
  = \mathbf{f}(t),
  \qquad t \in (0,t_{\mathrm{f}}],
  \label{eq:semi-discrete-ode}
\end{equation}
with initial condition
\begin{equation}
  \mathbf{u}(0) = \mathbf{u}_{0},
  \label{eq:semi-discrete-ic}
\end{equation}
where $\mathbf{u}(t) = (u_{\lambda}(t))_{\lambda \in \Lambda_{J}}$ is the
vector of expansion coefficients. 
The entries of the mass matrix $\mathbf{M}$, stiffness matrix $\mathbf{K}$
and load vector $\mathbf{f}(t)$ are given by
\begin{align}
  M_{\lambda\mu}
  &:= m(\psi_{\mu},\psi_{\lambda}),
  \label{eq:mass-matrix}
  \\
  K_{\lambda\mu}
  &:= a(\psi_{\mu},\psi_{\lambda}),
  \label{eq:stiffness-matrix}
  \\
  f_{\lambda}(t)
  &:= \ell(\psi_{\lambda},t),
  \label{eq:load-vector}
\end{align}
for all $\lambda,\mu \in \Lambda_{J}$. 
Because of the compact support of the basis functions and the locality of
the operators involved in \eqref{eq:bilinear-form}--\eqref{eq:mass-form},
each row of $\mathbf{M}$ and $\mathbf{K}$ contains only a small number of
non-zero entries. 
Moreover, it has been shown that the stiffness matrix, when expressed in a
wavelet basis, is not only sparse but also \emph{compressible} in the
sense that it can be well approximated by matrices with a prescribed
number of non-zero entries per row and column [13,14]. 
This property underlies the efficiency of wavelet--based solvers for
elliptic and parabolic problems.

In the time direction we discretise \eqref{eq:semi-discrete-ode} by an
implicit backward Euler scheme with time step $\Delta t > 0$. 
For $n=0,1,\dots,N_{t}-1$ and $t_{n} = n\,\Delta t$, this leads to the
recurrence
\begin{equation}
  \bigl(\mathbf{M} + \Delta t\,\mathbf{K}\bigr)\mathbf{u}^{n+1}
  = \mathbf{M}\,\mathbf{u}^{n} + \Delta t\,\mathbf{f}^{n+1},
  \label{eq:backward-euler}
\end{equation}
where $\mathbf{u}^{n} \approx \mathbf{u}(t_{n})$ and
$\mathbf{f}^{n+1} \approx \mathbf{f}(t_{n+1})$. 
The unconditional stability of backward Euler
[15] makes this combination particularly attractive in
situations where the spatial resolution is strongly refined near
interfaces or inclusions, leading to stiffness in the semi-discrete
system.

\subsection{Adaptive refinement driven by wavelet coefficients}

A major advantage of working with wavelet bases is the possibility of
constructing adaptive discretisations in a natural and inexpensive way.
Because wavelets are organised by scale and location, the magnitude of a
coefficient $u_{\lambda}(t)$ carries information about the local
regularity of the temperature field $T_{J}(\cdot,t)$ in the region where
$\psi_{\lambda}$ is supported. 
Small coefficients indicate locally smooth behaviour, whereas large
coefficients signal steep gradients or other fine-scale features [11,13].

We exploit this observation by introducing, at each time level, an active
index set $\Lambda_{J}^{\mathrm{act}}(t) \subset \Lambda_{J}$ and
representing the numerical solution by
\[
  T_{J}^{\mathrm{act}}(x,t)
  = \sum_{\lambda \in \Lambda_{J}^{\mathrm{act}}(t)}
      u_{\lambda}(t)\,\psi_{\lambda}(x).
\]
Starting from an initial set $\Lambda_{J}^{\mathrm{act}}(0)$ determined by
the expansion of $T_{0}$, the set is updated according to a simple
thresholding strategy: after each time step we mark all indices
$\lambda \in \Lambda_{J}$ such that
$|u_{\lambda}(t_{n+1})| \ge \varepsilon_{\mathrm{tol}}$ for a prescribed
tolerance $\varepsilon_{\mathrm{tol}} > 0$, together with a small
\emph{neighbourhood} of these indices in scale and position to maintain
stability of the scheme. The active set is updated as part of the definition for the next approximation space. The automatic refinements defined by this procedure develop the majority of resolution around material interfaces, inclusions, and the regions with rapidly developing transient boundary layers while leaving other relatively uniform portions of the domain in a relatively coarse level of resolution. The adaptive wavelet-Galerkin scheme can be viewed from an algorithmic perspective as a special case in the overall framework of adaptive wavelet-based techniques addressed in papers [11, 13, 16] for the solving of elliptic/parabolic mathematical problems. The aforementioned framework developed in these papers provides a very rigorous description of convergence and an upper and lower bound complexity for these types of adaptive wavelet schemes. The scheme exhibits a logarithmic growth in the number of degrees-of-freedom as a result of the use of wavelet decomposition, and when applied to composite materials, the scheme retains the ability to model complex microstructural details such as abrupt transitions across parallel layers or gradual transitions within multi-material composite structures using near-optimal levels of computational complexity. Numerical examples demonstrating these properties will be presented in the next section.
\section{Time Discretization and Algorithmic Workflow}
\label{sec:time-discretization}

The wavelet--Galerkin semi-discretisation introduced in
Section~\ref{sec:wavelet-galerkin} leads to the system of ordinary
differential equations
\begin{equation}
  \mathbf{M}\,\dot{\mathbf{u}}(t) + \mathbf{K}\,\mathbf{u}(t)
  = \mathbf{f}(t),
  \qquad t \in (0,t_{\mathrm{f}}],
  \label{eq:semi-discrete-ode-2}
\end{equation}
supplemented with the initial condition $\mathbf{u}(0) = \mathbf{u}_{0}$.
Here $\mathbf{M}$ and $\mathbf{K}$ are the mass and stiffness matrices
defined in \eqref{eq:mass-matrix}--\eqref{eq:stiffness-matrix}, while
$\mathbf{f}(t)$ collects the load contributions
\eqref{eq:load-vector}. 
In this section we describe the fully discrete time-stepping scheme and
summarise the resulting adaptive algorithm.

\subsection{Fully discrete backward Euler scheme}

To advance the semi-discrete solution in time we adopt a backward Euler
discretisation. 
For a given time step $\Delta t > 0$ we introduce the grid points
$t_{n} = n\,\Delta t$, $n = 0,1,\dots,N_{t}$, with $N_{t}\Delta t = t_{\mathrm{f}}$.
Let $\mathbf{u}^{n}$ denote the approximation to $\mathbf{u}(t_{n})$ and
$\mathbf{f}^{n+1}$ an approximation to the load vector at time
$t_{n+1}$. 
Applying backward Euler to \eqref{eq:semi-discrete-ode-2} yields
\begin{equation}
  \mathbf{M}\,\frac{\mathbf{u}^{n+1}-\mathbf{u}^{n}}{\Delta t}
  + \mathbf{K}\,\mathbf{u}^{n+1}
  = \mathbf{f}^{n+1},
  \qquad n = 0,1,\dots,N_{t}-1.
\end{equation}
Rearranging the terms we obtain the linear system
\begin{equation}
  \mathbf{A}\,\mathbf{u}^{n+1}
  = \mathbf{r}^{n},
  \qquad
  \mathbf{A} := \mathbf{M} + \Delta t\,\mathbf{K},
  \qquad
  \mathbf{r}^{n} := \mathbf{M}\,\mathbf{u}^{n} + \Delta t\,\mathbf{f}^{n+1}.
  \label{eq:backward-euler-fully-discrete}
\end{equation}
Under the assumptions of Section 2, the mass
matrix $\mathbf{M}$ is symmetric positive definite and the stiffness
matrix $\mathbf{K}$ is symmetric positive semi-definite. 
Consequently, for any $\Delta t > 0$ the matrix
$\mathbf{A} = \mathbf{M} + \Delta t\,\mathbf{K}$ remains symmetric
positive definite, and backward Euler is unconditionally stable in the
sense of energy estimates in the corresponding discrete norm [15,17]. 
This property is particularly relevant in the present setting, where the
adaptive wavelet refinement may lead to very fine local resolution near
interfaces and hence to a stiff system.

In the numerical experiments reported below we use a constant time step
$\Delta t$ for simplicity. 
More sophisticated strategies, in which $\Delta t$ is adapted in time
based on an error indicator or on physical time scales, could be
incorporated without altering the spatial discretisation [18]. 
Such extensions are left for future work.

\subsection{Solution of the linear systems}
At every time level $t_{n+1}$, we need to solve for $\mathbf{u}^{n+1}$ by solving equation (24) , which is a discretized form of the continuum PDE problem. The application of wavelet basis functions provides compact support for the wavelet coefficients; therefore, $M$ and $K$ will also be sparse, as will $A$. Furthermore, since the stiffness matrix expressed in a wavelet basis is known to be compressible, it is possible to obtain an accurate representation with matrices with uniformly bounded numbers of non-zero entries per row and column [13,14]. Because of this type of structure, the system is particularly well suited for an iterative solver approach. In our implementation, we use a preconditioned conjugate gradient (PCG) method for solving equation (24) , which is symmetric and positive definite. A preconditioner can be generated by using the diagonal of A or via a multilevel wavelet approximation using the wavelet scaling functions, leading to uniformly bounded condition numbers for the preconditioned operator at all levels of refinement [11,13]. As a result, our overall computational complexity remains close to linear with respect to the number of active degrees of freedom, regardless of both the local mesh sizes at sharp interfaces and the highest active resolution level. As an alternative, depending on the size and availability of the relevant libraries, one could also use either a sparse direct solver or a Krylov subspace method in the case of non-symmetric formulations [15]; however, for the category of composite heat-conduction problems studied here, the use of a symmetric wavelet-Galerkin discretisation in combination with PCG has been adequate.
\subsection{Algorithmic workflow with adaptive refinement}

We summarise here the overall algorithmic structure of the adaptive
wavelet--Galerkin scheme used in the numerical experiments. 
Let $\Lambda_{J}^{\mathrm{act}}(t_{n})$ denote the active index set at
time $t_{n}$ and $\varepsilon_{\mathrm{tol}}$ the threshold parameter
governing the refinement policy.

\begin{enumerate}
  \item \textbf{Initialisation.}  
        Choose the final time $t_{\mathrm{f}}$, the time step $\Delta t$,
        the maximal wavelet level $J$ and the threshold
        $\varepsilon_{\mathrm{tol}}$.  
        Compute the wavelet expansion of the initial condition $T_{0}$,
        determine the corresponding coefficient vector $\mathbf{u}^{0}$ and
        define the initial active set
        $\Lambda_{J}^{\mathrm{act}}(t_{0})$ as the collection of indices
        with $|u_{\lambda}^{0}| \ge \varepsilon_{\mathrm{tol}}$ plus a
        small neighbourhood in scale and position.

  \item \textbf{Assembly on the active set.}  
        Restrict the basis to the active indices
        $\Lambda_{J}^{\mathrm{act}}(t_{n})$ and assemble the associated
        submatrices $\mathbf{M}^{(n)}$ and $\mathbf{K}^{(n)}$, together
        with the load vector $\mathbf{f}^{(n+1)}$ at time $t_{n+1}$.  
        Form the system matrix
        $\mathbf{A}^{(n)} = \mathbf{M}^{(n)} + \Delta t\,\mathbf{K}^{(n)}$.

  \item \textbf{Time step.}  
        Solve
        \[
          \mathbf{A}^{(n)} \,\mathbf{u}^{n+1}
          = \mathbf{M}^{(n)}\,\mathbf{u}^{n}
            + \Delta t\,\mathbf{f}^{(n+1)}
        \]
        by the preconditioned conjugate gradient method up to a prescribed
        relative residual tolerance.

  \item \textbf{Coefficient-based adaptivity.}  
        Inspect the magnitudes of the updated coefficients
        $\{u_{\lambda}^{n+1}\}_{\lambda \in \Lambda_{J}^{\mathrm{act}}(t_{n})}$.
        Mark all indices with
        $|u_{\lambda}^{n+1}| \ge \varepsilon_{\mathrm{tol}}$ as
        \emph{essential} and enlarge this set by adding a small
        neighbourhood in scale and position to maintain stability of the
        representation.  
        This defines the next active index set
        $\Lambda_{J}^{\mathrm{act}}(t_{n+1})$.

  \item \textbf{Loop in time.}  
        If $t_{n+1} < t_{\mathrm{f}}$, set $n \leftarrow n+1$ and repeat
        steps 2--4.  
        Otherwise, stop and post-process the final temperature field and,
        if needed, derived quantities such as averaged heat fluxes or
        effective thermal conductivities.
\end{enumerate}
This approach ensures that locations where the temperature field contains fine-scale features (e.g. steep temperature changes due to material interfaces; prominent boundary layers, and locally concentrated heat) are the only areas where extra degrees of freedom will be added. In the examples presented, we see that the discretizations produced using this technique, on average, are much more economical than the uniform-resolution approaches, while maintaining the inherent advantages of the backward Euler time-stepping method including stability and robustness.
\vspace{1cm}
\section{Numerical Experiments and Performance Analysis}
\label{sec:numerical-experiments}
The goal of this section is to demonstrate how the adaptive wavelet-Galerkin method reacts to three examples of heterogeneous material behaviour due to heat conduction; I.e., the three types of problems presented shape beneath the surface: the first example is a laminated composite sheet; the second example is the presence of circular holes in an isotropic material; and the third example is a material with varying physical/thermal properties. These three represent typical configurations used in many thermal simulations of composite structures [2,9]. The objectives of this evaluation are: a) To determine if the adaptive wavelet-Galerkin method yields accurate solutions compared to existing reference solutions; and b) To estimate/quantify how much less memory and processing power is required with adaptively sampled solutions versus uniformly sampled solutions, based on the evaluation above.
The computational domain used for all calculations unless specified is the unit square $\Omega = (0,1)^{2}$, with a constant or smoothly varying conductivity tensor. The material parameterization does not include values for the material density $\rho$ and the heat capacity $c_{p}$; all of which have been set as unitless for all phases. Time is non-dimensionalized with respect to the characteristic diffusion time, which is the square of the length $L$ and divided by the value of the reference conductivity kref. The value of $L$  and $k_{\mathrm{ref}}$ are both established as characteristic parameters. The ambient temperature is set to zero, and Dirichlet or Neumann boundary conditions have been assigned; hence, the Robin term in equations (9) and (10) has been set to zero. In every test case, we compare the adaptive wavelet solution against a reference finite element approximation created from a significantly refined mesh and calculate the $L^{2}/k_{\mathrm{ref}}$ error associated with the temperature and the $H^{1}$ semi-norm of that error [4].
\vspace{1cm}
\subsection{Test case I: layered composite slab}
The first configuration is a two-phase layered composite, which provides a
clean benchmark for the resolution of sharp conductivity jumps across
planar interfaces. 
The domain $\Omega$ is partitioned into two horizontal layers,
\[
  \Omega_{1} = (0,1)\times(0,\tfrac{1}{2}), 
  \qquad
  \Omega_{2} = (0,1)\times(\tfrac{1}{2},1),
\]
with conductivity tensors
$\mathbf{K}_{1} = k_{1}\,\mathbf{I}$ in $\Omega_{1}$ and
$\mathbf{K}_{2} = k_{2}\,\mathbf{I}$ in $\Omega_{2}$, where
$k_{2}/k_{1} \gg 1$ models a high-conductivity layer embedded in a less
conductive matrix. 
The lower edge $y=0$ is held at a fixed temperature $T = 1$, the upper
edge $y=1$ is held at $T = 0$ and the lateral boundaries $x=0$ and $x=1$
are thermally insulated,
\[
  T(x,0,t) = 1,
  \quad
  T(x,1,t) = 0,
  \quad
  \mathbf{n}\cdot\mathbf{K}\nabla T = 0
  \ \text{on } x=0,1.
\]
\begin{figure}[H]
  \centering
  \includegraphics[width=0.45\textwidth]{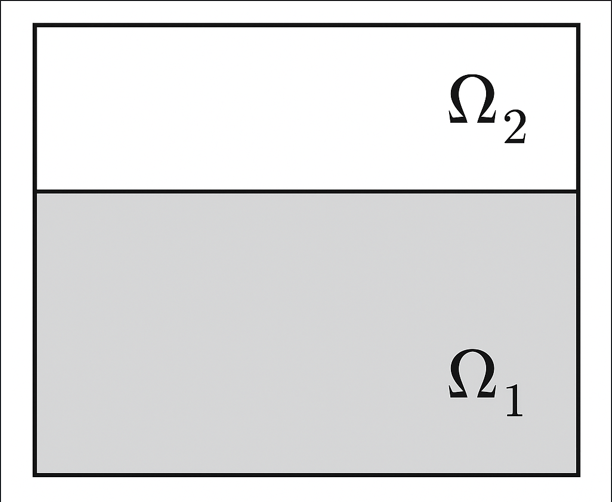}
  \vspace{1cm}
  \caption{Geometry of the two--phase layered composite slab used in Test case~I. 
  The domain $\Omega$ is split into a lower layer $\Omega_{1}$ and an upper layer 
  $\Omega_{2}$, separated by a planar interface at $y = 1/2$ with perfect thermal contact.}
  \label{fig:layered}
\end{figure}
The initial state is defined to be $T_{0}(x) \equiv 0$ in $\Omega$, which results in a transient thermal layer developing from the initial condition and relaxing to the steady-state solution. The discontinuity in the normal component of $\nabla T$ across the $y=\tfrac{1}{2}$ interface is also a result of the layered geometry, whereas the temperature value remains continuous across the interface. Using a uniform mesher will not give you the correct results, as you must be able to capture the discontinuity in the normal component of the gradient by placing mesh lines along the interface and using sufficiently small elements to represent this thermal layer. An adaptive wavelet scheme creates wavelet basis functions that automatically activate at the interface. The coefficients of the wavelets that are active at the interface will be refined, while the wavelets that are further away from the interface will have coarser representations. The active wavelet pattern will consist of a band of refined coefficients with a decreasing amount of refinement as the distance from the interface increases. For tolerances $\varepsilon_{\mathrm{tol}}$ that match the $L^{2}$-error of a fine uniform grid, the number of active degrees of freedom will be substantially fewer than the full dimensionality of $V_{J}$ , while still producing a thermal layer with no apparent spurious oscillations in the temperature profiles along vertical cuts.
\vspace{1cm}
\subsection{Test case II: circular inclusion in a matrix}

The second test case addresses a matrix-inclusion geometry with curved
interfaces, which is more challenging for mesh-based methods that require
either body-fitted grids or unfitted interface treatments. 
Here the unit square contains a circular inclusion of radius $r=0.2$
centred at $(0.5,0.5)$, denoted by $\Omega_{\mathrm{inc}}$, embedded in a
matrix region $\Omega_{\mathrm{m}} = \Omega \setminus \overline{\Omega_{\mathrm{inc}}}$. 
The conductivity is taken as $\mathbf{K} = k_{\mathrm{m}}\mathbf{I}$ in
the matrix and $\mathbf{K} = k_{\mathrm{inc}}\mathbf{I}$ inside the
inclusion, with a moderate contrast $k_{\mathrm{inc}}/k_{\mathrm{m}}$ to
avoid excessively sharp boundary layers.
\begin{figure}[H]
  \centering
  \includegraphics[width=0.45\textwidth]{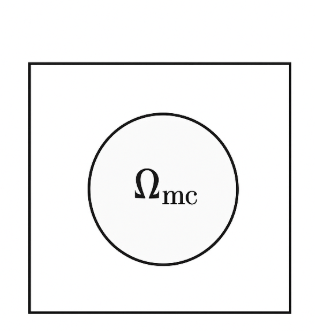}
  \caption{Geometry of the matrix--inclusion configuration used in Test case~II. 
  A circular inclusion $\Omega_{\mathrm{inc}}$ of radius $r$ is embedded in the matrix 
  domain $\Omega_{\mathrm{m}}$, with distinct thermal conductivities in the two regions.}
  \label{fig:inclusion}
\end{figure}

We impose a unit temperature on the left boundary and zero temperature on
the right boundary,
\[
  T(0,y,t) = 1,
  \qquad
  T(1,y,t) = 0,
\]
while the top and bottom edges are insulated,
\[
  \mathbf{n}\cdot\mathbf{K}\nabla T = 0
  \quad \text{on } y=0,1.
\]
The initial condition is set such that  $T_{0}(x) \equiv 0$. The resultant temperature profile has curved isotherms which are bent due to either more or less conductivity (depending on selection of parameters), resulting in an increase in heat flux near theicclusion. In this setup, the adaptive wavelet scheme creates a different annular enhancement pattern around the inclusion. While only a small pool of basis functions were active at intermediate and fine scales, the coarse scale accurately approximates the temperature distribution away from the inclusion. Comparison with the reference finite element solution has been made to determine that the wavelet-Galerkin solution is able to achieve an equal level of accuracy for global temperature distribution and localized heat flux perturbation near the intervening mass with a reduced number of active degrees of freedom. The $H^{1}$-seminorm for the error decreases consistently as $\varepsilon_{\mathrm{tol}}$  is reduced; whilst, convergence rates were also shown to correspond to the order of regularity for the actual solution underneath.
\subsection{Test case III: functionally graded medium}

The third example concerns a functionally graded medium (FGM), which
combines a smooth variation of material properties with localised
boundary layers. 
We consider a coating layer occupying the upper part of the domain,
\[
  \Omega_{\mathrm{fgm}} = (0,1)\times(\tfrac{1}{2},1),
\]
with a conductivity that varies continuously in the vertical direction
according to
\[
  \mathbf{K}(x,y) = k_{\mathrm{m}}\mathbf{I}
  \quad \text{for } 0 < y < \tfrac{1}{2},
  \qquad
  \mathbf{K}(x,y) = k_{\mathrm{m}}\Bigl(1 + \alpha (2y-1)\Bigr)\mathbf{I}
  \quad \text{for } \tfrac{1}{2} \le y \le 1,
\]
where $\alpha$ controls the grading intensity [17,18,19]. 
This example illustrates how the electrical conductivity of a metallic substrate can be affected by a ceramic layer whose conductivity either increases (or decreases) continuously throughout the thickness. The conditions at the boundaries of this example are identical to when the metallic substrate is layered, with the bottom of the ceramic layer held at a constant temperature $T=1$, the top of the ceramic layer kept at a lower temperature $T=0$ and both of the sides insulated.
\begin{figure}[H]
  \centering
  \includegraphics[width=0.45\textwidth]{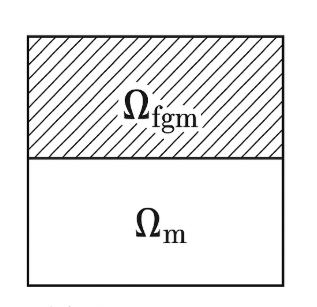}
  \caption{Geometry of the Case III functionally graded material. The upper layer of the FGM $\Omega_{\mathrm{fgm}}$ is a functionally graded coating that has a gradual change in thermal conductivity, while the lower section of the material $\Omega_{\mathrm{m}}$ is a homogeneous substrate to which it is bonded.}
  \label{fig:fgm}
\end{figure}
The conductivity gradient in the graded layer is instead continuous, resulting in a smooth temperature gradient without abrupt changes within the layer. However, since there exists a high thermal gradient close to the cooler boundary at $y=1$, there needs to be sufficient resolution in a thin boundary layer. As so, the resolution of the adaptive wavelet method reflects the geometry of this particular problem: the densest area is at the cooled area and gradually coarsens as you go further into the interior of the graded layer. The distribution of the active coefficient within the graded layer follows the spatial distribution of $\mathbf{K}(x,y)$. Thus, more coefficients will have more freedom of movement where effective diffusivity is changing most abruptly. The semi-analytical reference solution, which is produced by removing the lateral coefficient changes in the $K$ value, and was reduced to a one-dimensional case, supports the validity of the temperature profile determination from the adaptive wavelet method.

\subsection{Discussion of accuracy and efficiency}
The patterns observed in all three test cases for the adaptive wavelet–Galerkin solution method were similar. For a given tolerance level of $\varepsilon_{\mathrm{tol}}$, the combined $L^{2}$ and $H^{1}$ error of the adaptive solution is similar to that of the uniform wavelet discretisation at the finest resolution, $J$. However, the number of active wavelet coefficients required to represent the adaptive wavelet discretisation for $\varepsilon_{\mathrm{tol}}$ was significantly less than the maximum dimension of the uniform wavelet discretisation $V_{J}$. The reduction in the number of coefficients occurs most significantly in test cases 2 and 3, where in both examples the solution is comparatively uniform away from the much lower dimensional, relatively small set of interfaces. The reduction in degree of control is less marked for test case 1, where there is a smooth rate of change of temperature across the full thickness of the graded thermal coating. From an algorithmic perspective, the wavelet representation has both a sparsity and a compressibility feature which results in an assembly/sourcing phase linearly or closely approximating lines to that of the total number of variables being solved, as would be anticipated from the theoretical results presented for adaptive wavelet algorithms [11,13,16]. The preconditioned conjugate gradient scheme will converge in a few cycles that, for all practical reasons, is not dependent on the size of the finite element mesh. Therefore, it follows that the final computation time is primarily spent on the multiplication of M \& K by vectors. Based on the findings, it can be concluded that the adaptive wavelet–Galerkin methodology provides a comparable option to traditional finite-element approaches to thermal numerical modelling of composite materials where the distribution of materials, structures \& their geometrical forms yield a steep excess of thermal energy gradients.

\section{Conclusions and Future Work}
\label{sec:conclusions}

\subsection{Conclusions}

The adaptive wavelet-Galerkin technique has been developed through research carried out to assist in the solution of transient heat conduction in composite materials and functionally graded structures. The work involved expressing the anisotropic diffusion problem with mixed Dirichlet-Neumann boundary conditions in both strong and variational form before employing compactly supported wavelet basis functions for space discretisation, thereby creating a numerical representation of multi-layered and functionally graded structures. This semi-discrete model is unconditionally stable in respect of time when solved using a backward in time scheme; thus allowing for an efficient computational solution via preconditioned Krylov based solvers as the stiffness matrices associated with wavelets tend to be relatively sparse and contain multiple levels. Numerical experiment results demonstrate that the proposed method can achieve good accuracy with comparatively few DOF when compared to regular (standard) finite element and finite difference discretisation methods performed on uniform meshes. The ability to adaptively improve a solution's accuracy using a local wavelet threshold means that the solution can be improved in regions where there is a large change in gradient such as near the edge of material interfaces or where there is a steep temporal change in temperature. Additionally, it has been shown that an accurate approximation of a solution with a long time history and a realistic representation of the way heat dissipates is feasible, as demonstrated in the benchmarking of the proposed method with respect to dynamic loads and multiple material types. Beyond providing improved measures of error quantitatively, wavelet representations can also provide a qualitative advantage. Each subband of a specific wavelet decomposition can correspond to different spatial/temporal scales of the temperature field, therefore, it is possible to associate the effects of certain microstructure characteristics (layer thickness; gradation profile; inclusion geometry) on the macrothermal response of the system with the wavelet decomposition that is produced. The multiscale nature of this type of representation has the natural advantage of allowing for the development of reduced-order models for the purposes of designing the same thermal systems more efficiently.

\subsection{Future work}
There are many examples of further work that could be completed using the current framework. We will be working on the modelling aspect of the formulation to extend the wavelet–Galerkin approach to include fully 3D problems, strongly anisotropic conductivity tensors, and coupled thermomechanical responses. Temperature dependent material properties and thermal stresses will all be included in the generalised framework. Additionally, we can incorporate various nonlinear effects at the PDE level, but still employ the same adaptive wavelet technique for the spatial discretisation. The second major area of research will involve integrating the current solver more formally with the methods of homogenisation and effective property calculation. By linking the current solver to periodic or stochastic representative volume elements, high fidelity databases for effective thermal conductivity can be created from many different microstructural configurations for heterogeneous composites. These databases can be processed to identify scaling laws and create surrogate constitutive relationships for design purposes, similar to those obtained from current multiscale studies on woven and porous materials [17,18,19,20,21]. A promising direction for future research, and ultimately, the emphasis of our follow-up study will be the development of a hybrid wavelet and machine learning methodology to facilitate rapid in-situ prediction of the effective thermal properties of composites. The wavelet solver developed in this work will be used offline to create large datasets, each containing a set of materials and geometrical descriptors as well as a reduced number (but still informative) wavelet coefficients corresponding to the composite's temperature field and the keff value for the composite. Supervised machine learning models, such as Gaussian process regression, random forests and neural networks, can be trained using these compression features for the development of approximative prediction models for the keff. This approach will reduce or eliminate the time spent predicting keff using pure physical modelling as well as align with recent research in data-driven prediction of the thermal transport in heterogeneous materials. The hybrid wavelet and machine learning methodology may be employed for developing real-time design and optimisation tools for composite thermal systems while benefiting from the increased level of detail associated with the wavelet-Galerkin methodology.

\bibliographystyle{unsrtnat}

\end{document}